\documentclass{svmult}
 \usepackage{amsmath,amssymb,bbm}
 \smartqed

\def\PP{\mathbb P}
\def\PE{\mathbb{E}}

\def\Pdouble{\check{P}}
\def\PPdouble{\check{\PP}}
\def\PEdouble{\check{\PE}}

\def\Ptriple{\tilde{P}}
\def\PPtriple{\tilde{\PP}}
\def\PEtriple{\tilde{\PE}}

\newtheorem{step}{Step}

\begin{document}

\title*{Subgeometric ergodicity of Markov chains} \author{Randal
  Douc\inst{1} \and Eric Moulines\inst{2} \and Philippe
  Soulier\inst{3} }

\institute{CMAP, Ecole Polytechnique, 91128 Palaiseau Cedex, France
  \texttt{douc@cmap.polytechnique.fr} \and D\'epartement TSI, Ecole
    nationale sup\'erieure des T\'el\'ecommunications, 46 rue
    Barrault, 75013 Paris, France \texttt{moulines@tsi.enst.fr}
  \and Equipe MODAL'X, Universit\'e de Paris X Nanterre, 92000
  Nanterre, France \texttt{philippe.soulier@u-paris10.fr}}

  \maketitle

\section{Introduction}
Let $P$ be a Markov tranition kernel on a state space ${\mathsf{X}}$ equipped
with a countably generated $\sigma$-field $\mathcal{X}$.  For a control
function $f: {\mathsf{X}} \to [1,\infty)$, the $f$-\textit{total variation}
or $f$-\textit{norm} of a signed measure $\mu$ on $\mathcal{X}$ is defined
as
\[
\| \mu \|_f := \sup_{|g| \leq f} | \mu(g) | \; .
\]
When $f \equiv 1$, the $f$-norm is the total variation norm, which is
denoted $\| \mu \|_{\mathrm{TV}}$.  Assume that $P$ is aperiodic
positive Harris recurrent with stationary distribution $\pi$. Then the
iterated kernels $P^n(x,\cdot)$ converge to $\pi$. The rate of
convergence of $P^n(x,.)$ to $\pi$ does not depend on the starting
state $x$, but exact bounds may depend on $x$. Hence, it is of
interest to obtain non uniform or quantitative bounds of the following
form
\begin{gather}  \label{eq:protoresult}
  \sum_{n=1}^\infty r(n) \| P^n(x,\cdot) - \pi \|_f \leq g(x) \; , \quad
  \text{for all $x \in {\mathsf{X}}$}
\end{gather}
where $f$ is a control function, $\{ r(n) \}_{n \geq 0}$ is a non-decreasing
sequence, and $g$ is a nonnegative function which can be computed explicitly.

As emphasized in \cite[section 3.5]{roberts:rosenthal:2004}, quantitative
bounds have a substantial history in Markov chain theory. Applications are
numerous including convergence analysis of Markov Chain Monte Carlo (MCMC)
methods, transient analysis of queueing systems or storage models, etc.  With
few exception however, these quantitative bounds were derived under conditions
which imply geometric convergence, \textit{i.e.} $ r(n) = \beta^n$, for some
$\beta>1$ (see for instance \cite{meyn:tweedie:1994}, \cite{rosenthal:1995},
\cite{roberts:tweedie:1999}, \cite{roberts:rosenthal:2004}, and
\cite{baxendale:2005}).

Geometric convergence does not hold for many chains of practical
interest.  Hence it is necessary to derive bounds for chains which
converge to the stationary distribution at a rate $r$ which grows to
infinity slower than a geometric sequence. These sequences are called
subgeometric sequences and are defined in
\cite{nummelin:tuominen:1983} as non decreasing sequences $r$ such
that $\log r(n)/n \downarrow 0$ as $n \to \infty$.  These
sequences include among other examples the polynomial sequences
$r(n)=n^\gamma$ with $\gamma>0$ and subgeometric sequences $r(n)
\E^{cn^\delta}$ with $c>0$ and $ \delta \in(0,1)$.

The first general results proving subgeometric rates of convergence
were obtained by \cite{nummelin:tuominen:1983} and later extended by
\cite{tuominen:tweedie:1994}, but do not provide computable
expressions for the bound in the rhs of \eqref{eq:protoresult}.  A
direct route to quantitative bounds for subgeometric sequences has
been opened by \cite{veretennikov:1997,veretennikov:1999}, based on
coupling techniques.  Such techniques were later used in specific
contexts by many authors, among others, \cite{fort:moulines:2000}
\cite{jarner:roberts:2001} \cite{fort:2001}
\cite{fort:moulines:2003:SPA}.

The goal of this paper is to give a short and self contained proof of
general bounds for subgeometric rates of convergence, under practical
conditions. This is done in two steps. The first one is
Theorem~\ref{theo:mainresult} whose proof, based on coupling, provides
an intuitive understanding of the results of
\cite{nummelin:tuominen:1983} and \cite{tuominen:tweedie:1994}. The
second step is the use of a very general drift condition, recently
introduced in \cite{douc:fort:moulines:soulier:2004}. This condition
is recalled in Section~\ref{sec:drift} and the bounds it implied are
stated in Proposition~\ref{prop:drift}.

This paper complements the works
\cite{douc:fort:moulines:soulier:2004} and
\cite{douc:moulines:soulier:2005}, to which we refer for applications
of the present techniques to practical examples.

\section{Explicit bounds for the rate of convergence}
The only assumption for our main result is the existence of a small
set.
\begin{enumerate}
\item \label{assumption:smallset} There exist a set $C \in
  \mathcal{X}$, a constant $\epsilon > 0$ and a probability measure
  $\nu$ such that, for all $x \in C$, $P(x,\cdot) \geq \epsilon
  \nu(\cdot)$.
\end{enumerate}
For simplicity, only one-step minorisation is considered in this
paper.  Adaptations to $m$-step minorisation can be carried out as in
\cite{rosenthal:1995} (see also \cite{fort:2001} and
\cite{fort:moulines:2003:SPA}).

Let $\Pdouble$ be a Markov transition kernel on ${\mathsf{X}} \times
{\mathsf{X}}$ such that, for all $A \in \mathcal{X}$,
\begin{align} \label{eq:defPdouble} & \Pdouble(x,x', A \times
  {\mathsf{X}})= P(x,A) \mathbbm{1}_{(C \times C)^c}(x,x') + Q(x,A)
  \mathbbm{1}_{C \times C}(x,x')
  \\
  & \Pdouble(x,x', {\mathsf{X}} \times A)= P(x', A) \mathbbm{1}_{(C
    \times C)^c}(x,x') + Q(x',A) \mathbbm{1}_{C \times C}(x,x')
\end{align}
where $A^c$ denotes the complementary of the subset $A$ and $Q$ is the
so-called residual kernel defined, for $x \in C$ and $A \in \mathcal{X}$
by
\begin{gather}
\label{eq:DefinitionResidualKernel}
Q(x,A) = \begin{cases}
(1-\epsilon)^{-1} \left( P(x,A) - \epsilon \nu(A) \right) & 0 < \epsilon < 1 \\
\nu(A) & \epsilon = 1
\end{cases}
\end{gather}
One may for example set
\begin{multline} \label{eq:bivariatekernel-independentcomponent}
  \Pdouble(x,x'; A \times A') = \\
  P(x,A) P(x',A') \mathbbm{1}_{(C \times C)^c}(x,x') + Q(x,A) Q (x',A)
  \mathbbm{1}_{C \times C}(x,x') \; ,
\end{multline}
but this choice is not always the most suitable; cf.
Section~\ref{sec:monotone}.  For $(x,x') \in {\mathsf{X}} \times {\mathsf{X}}$,
denote by $\PPdouble_{x,x'}$ and $\PEdouble_{x,x'}$ the law and the
expectation of a Markov chain with initial distribution $\delta_x
\otimes \delta_{x'}$ and transition kernel $\Pdouble$.

\begin{theorem} \label{theo:mainresult}
  Assume \ref{assumption:smallset}.  

\medskip 

\noindent
For any sequence $r\in\Lambda$, $\delta >0$ and all $(x,x') \in
{\mathsf{X}} \times {\mathsf{X}}$,
\begin{gather} \label{eq:sumtv}
  \sum_{n=1}^\infty r(n) \| P^n(x,\cdot) - P^n(x',\cdot)
  \|_{\mathrm{TV}} \leq (1+\delta) \PEdouble_{x,x'} \left[
    \sum_{k=0}^{\sigma} r(k) \right] + \frac{1-\epsilon}\epsilon M \; ,
\end{gather}
with $M = (1+\delta) \sup_{n\geq0} \left\{ R^* r(n-1) -
  \epsilon(1-\epsilon) \delta R(n) /(1+\delta) \right\}_+$ and $R^*=
\sup_{(y,y') \in {C \times C}} \PEdouble_{y,y'} \left[
  \sum_{k=1}^{\tau} r(k) \right]$.

\medskip 
\noindent
Let $W:{\mathsf{X}}\times{\mathsf{X}} \to [1,\infty)$ and $f$ be a non-negative
function $f$ such that $f(x) + f(x') \leq W(x,x')$ for all $(x,x') \in
{\mathsf{X}} \times {\mathsf{X}}$. Then,
\begin{gather} \label{eq:fnorm}
  \sum_{n=1}^\infty \| P^n(x,\cdot) - P^n(x',\cdot) \|_{f} \leq
  \PEdouble_{x,x'} \left[ \sum_{k=0}^{\sigma} W(X_k,X'_k)\right] +
  \frac{1-\epsilon}\epsilon W^* \; .
\end{gather}
with $W^* = \sup_{(y,y') \in {C \times C}} \PEdouble_{y,y'} \left[
  \sum_{k=1}^{\tau} W(X_k,X'_k) \right] $.
\end{theorem}

\begin{remark}
  Integrating these bounds with respect to $\pi(\D x')$ yields similar bounds for
  $\|P^n(x, \cdot) - \pi\|_{\mathrm{TV}}$ and $\|P^n(x, \cdot) - \pi\|_{f}$.
\end{remark}

\begin{remark}
  The trade off between the size of the coupling set and the
  constant~$\epsilon$ appears clearly: if the small set is big, then the chain
  returns more often to the small set and the moments of the hitting times can
  expected to be smaller, but the constant $\epsilon$ will be smaller. This
  trade-off is illustrated numerically in \cite[Section
  3]{douc:moulines:soulier:2005}.
\end{remark}

By interpolation, intermediate rates of convergence can be obtained.
Let $\alpha$ and $\beta$ be positive and increasing functions such
that, for some $0 \leq \rho \leq 1$,
\begin{gather} \label{eq:additivity}
  \alpha(u) \beta(v) \leq \rho u+ (1-\rho) v \; , \quad \text{for all
    $(u,v) \in {\mathbb{R}}^+ \times {\mathbb{R}}^+$} \; .
\end{gather}
Functions satisfying this condition can be obtained from Young's
inequality. Let $\psi$ be a real valued, continuous, strictly
increasing function on ${\mathbb{R}}^+$ such that $\psi(0)= 0$; then for all
$a,b > 0$,
$$
ab \leq \Psi(a) + \Phi(b) \; , \text{where} \quad
\Psi(a) = \int_0^a \psi(x) dx \quad \text{and} \quad
\Phi(b) = \int_0^b \psi^{-1}(x) dx \; ,
$$
where $\psi^{-1}$ is the inverse function of $\psi$. If we set
$\alpha(u) = \Psi^{-1}(\rho u)$ and $\beta(v)= \Phi^{-1}((1-\rho)v)$,
then the pair $(\alpha,\beta)$ satisfies \eqref{eq:additivity}. A
trivial example is obtained by taking $\psi(x)= x^{p-1}$ for some $p
\geq 1$, which yields $\alpha(u) = (p \rho u)^{1/p}$ and $\beta(u) =
(p (1- \rho) u/(p-1))^{(p-1)/p}$. Other examples are given in
Section~\ref{sec:drift}.
\begin{corollary} \label{coro:interpolation}
  Let $\alpha$ and $\beta$ be two positive functions satisfying
  \eqref{eq:additivity} for some $0 \leq \rho \leq 1$.  Then, for any
  non-negative function $f$ such that $f(x) + f(x') \leq \beta \circ
  W(x,x')$ and $\delta>0$, for all $x,x'\in{\mathsf{X}}$ and $n \geq 1$,
  \begin{multline} \label{eq:interpolationfnorm}
    \sum_{n=1}^\infty \alpha(r(n)) \| P^n(x,\cdot) - P^n(x',\cdot)
    \|_{f} \leq \rho (1+\delta) \PEdouble_{x,x'} \left[
      \sum_{k=0}^{\sigma} r(k) \right] \\
    + (1-\rho) \PEdouble_{x,x'} \left[ \sum_{k=0}^{\sigma} W(X_k,X'_k)
    \right] \frac{1-\epsilon}\epsilon \{\rho M + (1-\rho) W^*\} \; .
  \end{multline}
\end{corollary}

\subsection{Drift Conditions for subgeometric ergodicity}
\label{sec:drift}
The bounds obtained in Theorem~\ref{theo:mainresult} and
Corollary~\ref{coro:interpolation} are meaningful only if they are finite.
Sufficient conditions are given in this section in the form of drift
conditions.
The most well known drift condition is the so-called Foster-Lyapounov drift
condition which not only implies but is actually equivalent to geometric
convergence to the stationary distribution, cf. \cite[Chapter 16]{meyn:tweedie:1993}.
\cite{jarner:roberts:2001}, simplifying and generalizing an argument in
\cite{fort:moulines:2000}, introduced a drift condition which implies
polynomial rates of convergence. We consider here the following drift
condition, introduced in \cite{douc:fort:moulines:soulier:2004}, which allows to bridge
the gap between polynomial and geometric rates of convergence.

\noindent \textbf{Condition D}($\phi,V,C$): There exist a
function $V: {\mathsf{X}} \to [1,\infty]$, a concave monotone non decreasing
differentiable function $\phi: [1, \infty] \mapsto (0, \infty]$, a
measurable set $C$ and a constant $b>0$ such that
\begin{gather*}
  P V + \phi \circ V \leq V + b \mathbbm{1}_C.
\end{gather*}

If the function $\phi$ is concave, non decreasing and differentiable,
define
\begin{gather} \label{eq:defhphi}
  H_\phi(v) := \int_1^v \frac{dx}{\phi(x)}.
\end{gather}
Then $H_\phi$ is a non decreasing concave differentiable function on
$[1, \infty)$. Moreover, since $\phi$ is concave, $\phi'$ is non
increasing.  Hence $\phi(v) \leq \phi(1) + \phi'(1) (v-1)$ for all
$v\geq1$, which implies that $H_\phi$ increases to infinity. We can
thus define its inverse $H_\phi^{-1}: [0, \infty) \to [1, \infty)$,
which is also an increasing and differentiable function, with
derivative $(H_\phi^{-1})'(x) = \phi\circ H_\phi^{-1}(x)$. For
$k\in{\mathbb{N}}$, $z\geq0$ and $v\geq1$, define
\begin{align}
  &  r_\phi(z) := (H_\phi^{-1})'(z)  = \phi \circ H_\phi^{-1}(z) \; .
  \label{eq:ratefunction}%
\end{align}
It is readily checked that if $\lim_{ t \to \infty } \phi'(t)= 0$,
then $r_\phi \in \Lambda$, cf \cite[Lemma
2.3]{douc:fort:moulines:soulier:2004}.

Proposition~2.2 and Theorem~2.3 in \cite{douc:moulines:soulier:2005}
show that the drift condition ${\bf D}(\phi,V,C)$ implies that the
bounds of Theorem~\ref{theo:mainresult} are finite. We gather here
these results.
\begin{proposition} \label{prop:drift}
  Assume that Condition ${\bf D}(\phi,V,C)$ holds for some small set
  $C$ and that $\inf_{x\notin C} \phi\circ V(x)>b$. Fix some arbitrary
  $\lambda \in (0, 1 - b / \inf_{x\notin C} \phi\circ V(x))$ and
  define $W(x,x')= \lambda\phi(V(x) + V(x') -1)$. Define also $V^* =
  (1-\epsilon)^{-1} \sup_{y\in C} \left \{ P V(y) - \epsilon
    \nu(V)\right\}$. Let $\sigma$ be the hitting time of the set
  $C\times C$. Then
\begin{gather*}
  \PEdouble_{x,x'} \left[ \sum_{k=0}^{\sigma} r_\phi(k) \right] \leq 1
  + \frac{r_\phi(1)}{\phi(1)} \left\{ V(x) + V(x') \right\} \mathbbm{1}_{(x,x')
    \notin C \times C} \; ,
  \\
  \PEdouble_{x,x'} \left[ \sum_{k=0}^{\sigma} W(X_k,X'_k) \right] \leq
  \sup_{{(y,y') \in C \times C}} W(y,y') +
  \{V(x)+V(x')\} \mathbbm{1}_{(x,x') \notin C \times C} \; , \\
  R^* \leq 1 + \frac{r_\phi(1)}{\phi(1)} \left\{ 2 V^* - 1\right\}
  \\
  W^* \leq \sup_{{(y,y') \in C \times C}} W(y,y') + 2 V^* - 1 \; .
\end{gather*}
\end{proposition}

\begin{remark}
  The condition $\inf_{y \notin C} \phi \circ V(y) >b$ may not be
  fulfilled. If level sets $\{ V \leq d \}$ are small, then the set
  $C$ can be enlarged so that this condition holds.  This additional
  condition may appear rather strong, but can be weakened by using
  small sets associated to some iterate $P^m$ of the kernel (see
  \textit{e.g.}  \cite{rosenthal:1995}, \cite{fort:2001} and
  \cite{fort:moulines:2003:SPA}). 
\end{remark}

We now give examples of rates that can be obtained
by~(\ref{eq:ratefunction}).
\paragraph{Polynomial rates}
Polynomial rates of convergence are obtained when Condition ${\bf
  D}(\phi,V,C)$ holds with $\phi(v) = c v^\alpha$ for some $\alpha\in
[0,1)$ and $c \in(0,1]$.  The rate of convergence in total variation
distance is $ r_\phi(n) \propto n^{\alpha/(1-\alpha)}$ and the pairs
$(r,f)$ for which~(\ref{eq:interpolationfnorm}) holds are of the form
$(n^{ (1-p) \alpha /(1-\alpha)}, V^{\alpha p})$ for $p\in[0,1]$, or in
other terms, $(n^{\kappa -1}, V^{1-\kappa (1-\alpha)})$ for $1 \leq
\kappa \leq 1/(1-\alpha)$, which is Theorem 3.6 of
\cite{jarner:roberts:2001}.

It is possible to extend this result by using more general
interpolation functions. For instance, choosing for $b>0$, $\alpha(x)
= (1\vee\log(x))^b$ and $\beta(x) = x (1\vee\log(x))^{-b}$ yields the
pairs $(n^{ (1-p) \alpha/(1-\alpha)} \log^b (n), V^{\alpha p} (1+\log
V)^{-b})$, for $p\in[0,1]$.

\paragraph{Logarithmic rates  of convergence}
Rates of convergence slower than any polynomial can be obtained when
condition {\bf D}($\phi,V,C)$ holds with a function $\phi$ that
increases to infinity slower than polynomially, for instance $\phi(v)
= c(1+\log(v))^\alpha$ for some $\alpha \geq 0$ and $c\in(0,1]$. A
straightforward calculation shows that
\[
r_\phi(n) \asymp \log^\alpha(n) \; .
\]
Pairs for which~(\ref{eq:interpolationfnorm}) holds are thus of the
form $((1+\log(n))^{(1-p) \alpha}, (1+ \log(V))^{p\alpha})$.

\paragraph{Subexponential  rates of convergence}
Subexponential rates of convergence faster than any polynomial are
obtained when the condition {\bf D}($\phi,V,C$) holds with $\phi$ such
that $v/\phi(v)$ goes to infinity slower than polynomially.  Assume
for instance that $\phi$ is concave and differentiable on
$[1,+\infty)$ and that for large $v$, $\phi(v) = c v /\log^\alpha(v)$
for some $\alpha>0$ and $c>0$.  A simple calculation yields
\[
r_\phi(n) \asymp n^{-\alpha/(1+\alpha)} \exp\left( \{c(1+\alpha)
  n\}^{1/(1+\alpha)} \right) \; .
\]
Choosing $\alpha(x) = x^{1-p} (1\vee \log(x))^{-b}$ and $\beta(x)= x^p
(1\vee \log(x))^b$ for $p\in(0,1)$ and $b\in{\mathbb{R}}$; or $p=0$ and
$b>0$; or $p=1$ and $b<-\alpha$ yields the pairs
\[
n^{-(\alpha+b)/(1+\alpha)} \exp\left( (1-p)\{ c(1+\alpha)
  n\}^{1/(1+\alpha)}\right) \, , \; V^{p} (1+ \log V)^{b} \; .
\]

\subsection{Stochastically monotone chains}
\label{sec:monotone}
Let ${\mathsf{X}}$ be a totally ordered set and let the order relation be
denoted by $\preceq$ and for $a\in{\mathsf{X}}$, let $(-\infty,a]$ denote the
set of all $x\in{\mathsf{X}}$ such that $x\preceq a$. A transition kernel on
${\mathsf{X}}$ is said to be stochastically monotone if $x\preceq y$ implies
$P(x , (-\infty,a]) \geq P(y,(-\infty,a])$ for all $a\in{\mathsf{X}}$. If
Assumption~\ref{assumption:smallset} holds, for a small set
$C=(-\infty,a_0]$, then instead of defining the kernel $\Pdouble$ as
in~(\ref{eq:bivariatekernel-independentcomponent}), it is convenient
to define it, for $x,x'\in{\mathsf{X}}$ and
$A\in\mathcal{X}\otimes\mathcal{X}$, by
\begin{multline*}
  \Pdouble(x,x';A) = \mathbbm{1}_{(x,x') \notin C \times C} \int_0^1
  \mathbbm{1}_{A}(P^\leftarrow(x,u),P^\leftarrow(x',u)) \, du
  \\
  + \mathbbm{1}_{C \times C} (x,x') \; \int_0^1
  \mathbbm{1}_{A}(Q^\leftarrow(x,u),Q^\leftarrow(x',u)) \, du \; ,
\end{multline*}
where, for any transition kernel $K$ on ${\mathsf{X}}$,
$K^\leftarrow(x,\cdot)$ is the quantile function of the probability
measure $K(x,\cdot)$, and $Q$ is the residual kernel defined
in~(\ref{eq:DefinitionResidualKernel}). This construction makes the
set $\{(x,x')\in{\mathsf{X}}\times {\mathsf{X}}:\, x \preceq x'\}$ absorbing
for~$\Pdouble$. This means that if the chain $(X_n,X'_n)$ starts at
$(x_0,x'_0)$ with $x_0 \preceq x'_0$, then almost surely, $X_n \preceq
X'_n$ for all $n$. Let now $\sigma_C$ and $\sigma_{C\times C}$ denote
the hitting times of the sets $C$ and $C\times C$, respectively. Then,
we have the following very simple relations between the moments of the
hitting times of the one dimensional chain and that of the
bidimensional chain with transition kernel $\Pdouble$.  For any
sequence $r$ and any non negative function $V$ all $x\preceq x'$
\begin{gather*}
  \PEdouble_{x,x'} \left[ \sum_{k=0}^{\sigma_{C\times C}} r(k)
    V(X_k,X_k') \right] \leq \PE_{x'} \left[ \sum_{k=0}^{\sigma_C}
    r(k) V(X_k') \right] \;   .
\end{gather*}
A similar bound obviously holds for the return times. Thus, there only
remain to obtain bounds for this quantities, which is very
straightforward if moreover condition $\bf{D}(\phi,V,C)$ holds.
Examples of stochastically monotone chains with applications to
queuing and Monte-Carlo simulation that satisfy
condition~$\bf{D}(\phi,V,C)$ are given in \cite[section
3]{douc:moulines:soulier:2005}.

\section{Proof of Theorem~\ref{theo:mainresult}}

Define a transition kernel $\Ptriple$ on the space $\tilde X = {\mathsf{X}}
\times {\mathsf{X}} \times \{ 0, 1\}$ endowed with the product $\sigma$-field
$\tilde{\mathcal{X}}$, for any $x,x' \in {\mathsf{X}}$ and $A,A'\in
\mathcal{X}$, by
\begin{align}
  \Ptriple \left((x,x',0), A \times A' \times \{0\}\right) & = \{1 -
  \epsilon \mathbbm{1}_{C \times C}(x,x')\} \Pdouble((x,x'), A \times A')
  \;, \label{eq:TildeQ1}\\
  \Ptriple \left((x,x',0), A \times A' \times \{1\}\right) & =
  \epsilon \mathbbm{1}_{C \times C}(x,x') \nu_{x,x'}(A \cap A') \;,
  \label{eq:TildeQ2}
  \\
  \Ptriple\left((x,x',1), A \times A' \times \{1\}\right) & = P(x,A
  \cap A') \;.
  \label{eq:TildeQ3}
\end{align}
For any probability measure $\tilde \mu$ on $(\tilde{{\mathsf{X}}}, \tilde
{\mathcal{X}})$, let $\PPtriple_{\tilde\mu}$ be the probability
measure on the canonical space $(\tilde{{\mathsf{X}}}^{{\mathbb{N}}},
\tilde{\mathcal{X}}^{\otimes {\mathbb{N}}})$ such that the coordinate process
$\{\tilde X_k\}$ is a Markov chain with transition kernel $\tilde P$
and initial distribution $\tilde \mu$.  The corresponding expectation
operator is denoted by ${\PEtriple}_{\tilde{\mu}}$.

The transition kernel $\Ptriple$ can be described algorithmically.
Given $\tilde X_0=(X_0,X'_0,d_0)=(x,x',d)$, $\tilde
X_1=(X_{1},X'_{1},d_{1})$ is obtained as follows.
\begin{itemize}
\item If $d = 1$ then draw $X_{1}$ from $P(x,\cdot)$ and set
  $X'_{1}=X_{1}$, $d_{1} = 1$.
\item If $d=0$ and $(x,x') \in {C \times C}$, flip a coin with
  probability of heads $\epsilon$.
\begin{itemize}
\item[--] If the coin comes up heads, draw $X_{1}$ from $ \nu_{x,x'}$
  and set $X'_{1}=X_{1}$ and $d_{1}=1$.
\item[--] If the coin comes up tails, draw $(X_{1},X'_{1})$ from
  $\Pdouble(x,x'; \cdot)$ and set $d_{1}=0$.
\end{itemize}
\item If $d=0$ and $(x,x') \not\in {C \times C}$, draw
  $(X_{1},X'_{1})$ from $\Pdouble(x,x'; \cdot)$ and set $d_{1}= 0$.
\end{itemize}
The variable $d_n$ is called the \emph{bell variable}; it indicates
whether coupling has occurred by time $n$ ($d_n=1$) or not ($d_n=0$).
The first index $n$ at which $d_n=1$ is the coupling time;
\begin{gather*}
T = \inf \{ k \geq 1 : d_k = 1 \}.
\end{gather*}
If $d_n=1$ then $X_k = X'_k$ for all $k \geq n$.  This coupling
construction is carried out in such a way that under $\tilde \PP_{\xi
  \otimes \xi' \otimes \delta_0}$, $\{X_k\}$ and $\{ X'_k \}$ are
Markov chains with transition kernel $P$ with initial distributions
$\xi$ and $\xi'$ respectively.

\medskip
The main tool of the proof is the following relation between
$\PEtriple_{x,x',0}$ and $\PEdouble_{x,x'}$, proved in \cite[Lemma
1]{douc:moulines:rosenthal:2004}.  For any non-negative adapted
process $(\chi_k)_{k\geq 0}$ and $(x,x') \in {\mathsf{X}} \times {\mathsf{X}}$,
\begin{gather} \label{eq:stop}
  \tilde \PE_{x,x',0} [ \chi_{n} \mathbbm{1}_{ \{T > n \} } ] = \PEdouble_{x,x'}
  \left[\chi_{n} \, (1-\epsilon)^{N_{n-1}} \right] \; ,
\end{gather}
where $N_n = \sum_{i=0}^n \mathbbm{1}_{C \times C}(X_i,X'_i)$ is the number of
visits to ${C \times C}$ before time~$n$.

We now proceed with the proof of Theorem~\ref{theo:mainresult}.

\begin{step} Lindvall's inequality \cite{lindvall:1979,lindvall:1992}
  \begin{gather} \label{eq:borneseriefnorm}
    \sum_{k=0}^\infty r(k) \| P^k(x,\cdot) - P^k(x',\cdot) \|_f \leq
    \PEtriple_{x,x',0} \left [ \sum_{j=0}^{T-1} r(j) \; \{f(X_j)+
      f(X'_j)\} \right] \; .
\end{gather}
\end{step}
\begin{proof}
  For any measurable function $\phi$ such that $|\phi|\leq f$, and for
  any $(x,x') \in {\mathsf{X}}\times{\mathsf{X}}$ it holds that
\begin{align*}
  |P^k\phi(x)-P^k\phi(x')| & = \left| \PEtriple_{x,x',0}[\{\phi(X_k)
    - \phi(X_k')\} \mathbbm{1}_{\{d_k=0\}}] \right| \\
  & \leq \PEtriple_{x,x',0}[ \{f(X_k)+ f (X'_k)\} \mathbbm{1}_{\{T>k\}}] \; .
  \end{align*}
  Hence $\|P^k(x,\cdots) - P^k(x',\cdot)\|_f \leq \PEtriple_{x,x',0}[
  \{f(X_k)+ f (X'_k)\} \mathbbm{1}_{\{T>k\}}]$.  Summing over $k$ yields
  (\ref{eq:borneseriefnorm}). \qed
\end{proof}

\begin{step}
  Denote $W_{r,f} (x,x') = \PEdouble_{x,x'} \left[\sum_{k=0}^{\sigma}
    r(k) f(X_k,X'_k) \right]$ and $W^*(r,f) = \sup_{(x,x') \in C
    \times C} \left[\sum_{k=1}^{\tau} r(k) f(X_k,X'_k)
  \right]/r(0)$. Then
  \begin{multline} \label{eq:relfund}
    \PEtriple_{x,x',0} \left [ \sum_{k=0}^{T-1} r(k) f(X_k,X'_k)
    \right] \\
    \leq W_{r,f}(x,x') + \epsilon^{-1} (1-\epsilon) \, W^*_{r,f} \;
    \PEtriple_{x,x',0} [r(T-1)] \; .
\end{multline}
\end{step}
\begin{proof}
  Applying~(\ref{eq:stop}), we obtain
  \begin{align*}
    \PEtriple_{x,x',0} & \left [ \sum_{k=0}^{T-1} r(k) f(X_k,X'_k)
    \right]  = \sum_{k=0}^\infty \PEtriple_{x,x',0} \left [ r(k)
      f(X_k,X'_k) \mathbbm{1}_{\{T>k\}} \right] \\
    & = \sum_{k=0}^\infty \PEdouble_{x,x'} \left [ r(k) f(X_k,X'_k)
      (1-\epsilon)^{N_{k-1}} \right] \\
    & = \sum_{j=0}^\infty \sum_{k=0}^\infty (1-\epsilon)^j
    \PEdouble_{x,x'} \left [ r(k) f(X_k,X'_k) \mathbbm{1}_{\{N_{k-1}=j\}}
    \right]    \\
    & = W_{r,f}(x,x') + \sum_{j=1}^\infty \sum_{k=0}^\infty
    (1-\epsilon)^j \PEdouble_{x,x'} \left [ r(k) f(X_k,X'_k)
      \mathbbm{1}_{\{N_{k-1}=j\}} \right]
  \end{align*}
  For $j\geq0$, let $\sigma_j$ denote the $(j+1)$-th visit to $C
  \times C$. Then $N_{k-1} = j$ iff $\sigma_{j-1} < k \leq \sigma_j$.
  Since $r$ is a subgeometric sequence, $r(n+m)\leq r(n)r(m)/r(0)$,
  thus
  \begin{align*}
    \sum_{k=0}^\infty r(k) f(X_k,X'_k) \mathbbm{1}_{\{N_{k-1}=j\}} & =
    \sum_{k=\sigma_{j-1}+1}^{\sigma_j} r(k) f(X_k,X'_k) \\
    & = \sum_{k=1}^{\tau \circ \theta^{\sigma_{j-1}}}
    r(\sigma_{j-1}+k)
    f(X_k,X'_k)  \\
    & \leq \frac{r(\sigma_{j-1})}{r(0)} \left( \sum_{k=1}^{\tau \circ
        \theta^{\sigma_{j-1}}} r(k) f(X_k,X'_k) \right) \circ
    \theta^{\sigma_{j-1}} \; .
  \end{align*}
  Applying the strong Markov property yields
\begin{multline*}
  \PEtriple_{x,x',0} \left [ \sum_{k=0}^{T-1} r(k) f(X_k,X'_k)
  \right] \leq W_{r,f}(x,x') \\
  + (1-\epsilon) W^*(f,g) \sum_{j=0}^\infty (1-\epsilon)^j
  \PEdouble_{x,x'}[r(\sigma_j)] \; .
\end{multline*}
By similar calculations,~(\ref{eq:stop}) yields
\begin{align*}
  \PEtriple[r(T-1)] = \epsilon \sum_{j=0}^\infty (1-\epsilon)^j
  \PEdouble[r(\sigma_j)] \; ,
\end{align*}
which concludes the proof of~(\ref{eq:relfund}). \qed
\end{proof}
\begin{step} \emph{
  Applying~(\ref{eq:relfund}) with $r\equiv1$ yields~(\ref{eq:fnorm}).}
\end{step}
\begin{step} \emph{  If $r \in \Lambda$, then $\lim_{n \rightarrow \infty} r(n) /
    R(k) = 0$, with $R(0)=1$ and $R(n) = \sum_{k=0}^{n-1} r(k)$,
    $n\geq1$.  Thus we can define, for $r \in \Lambda$ and $\delta >
    0$}
\begin{gather*}
  M_\delta = (1+\delta) \sup_{n\geq0} \left\{ \epsilon^{-1}
    (1-\epsilon) W_{r,1}^* r(n-1) - \delta R(n) /(1+\delta) \right\}_+
  \; .
\end{gather*}
\end{step}
$M_\delta$ is finite for all $\delta>0$. This yields
\begin{align*}
  \PEtriple_{x,x',0} [R(T)] \leq (1+\delta) W_{r,1}(x,x') + M_\delta \; .
\end{align*}
Applying this bound with~(\ref{eq:borneseriefnorm})
yields~(\ref{eq:sumtv}).   \qed

\end{document}